\documentclass{elsart}
\usepackage{graphicx}
\usepackage{amssymb,amsfonts,bm}
\newcommand{\trace}{\mathop{\rm Tr}\nolimits}
\newcommand{\var}{\mathop{\rm Var}\nolimits}

\newcommand{\diag}{\mathop{\rm Diag}\nolimits}

\newcommand{\twomat}[4]{\left(\begin{array}{cc}{#1}&{#2}\\{#3}&{#4}\end{array}\right)}

\newcommand{\cX}{{\mathcal X}} 

\newcommand{\C}{{\mathbb{C}}}
\newcommand{\R}{{\mathbb{R}}}
\newcommand{\M}{{\mathbb{M}}}
\DeclareRobustCommand\openone{\leavevmode\hbox{\small1\normalsize\kern-.33em1}}
\newcommand{\id}{\mathrm{\openone}}

\newcommand{\be}{\begin{equation}}
\newcommand{\ee}{\end{equation}}
\newcommand{\bea}{\begin{eqnarray}}
\newcommand{\eea}{\end{eqnarray}}
\newcommand{\beas}{\begin{eqnarray*}}
\newcommand{\eeas}{\end{eqnarray*}}
\newtheorem{definition}{Definition}
\newtheorem{theorem}{Theorem}
\newtheorem{lemma}{Lemma}
\newtheorem{corollary}{Corollary}

\newtheorem{conjecture}{Conjecture}
\newtheorem{problem}{Problem}
\newcount\minute
\newcount\hour
\def\currenttime{%
    \minute\time
    \hour\minute
    \divide\hour60
    \the\hour:\multiply\hour60\advance\minute-\hour\the\minute}
\begin{document}
\begin{frontmatter}
\title{Variance bounds, with an application to norm bounds for commutators}
\author{Koenraad M.R. Audenaert}
\address{
Mathematics Department\\
Royal Holloway, University of London\\
Egham TW20 0EX, United Kingdom
}
\ead{koenraad.audenaert@rhul.ac.uk}
\date{\today, \currenttime}
\begin{abstract}
Murthy and Sethi (Sankhya Ser B \textbf{27}, 201--210 (1965)) gave a sharp upper bound
on the variance of a real random variable in terms of
the range of values of that variable. We generalise this bound to
the complex case and, more importantly, to the matrix case. In doing so,
we make contact with several geometrical and matrix analytical
concepts, such as the numerical range, and introduce the new concept of radius of a matrix.

We also give a new and
simplified proof for a sharp upper bound
on the Frobenius norm of commutators recently proven
by B{\"o}ttcher and Wenzel
(Lin.\ Alg.\ Appl.\ \textbf{429} (2008) 1864--1885) and point out that
at the heart of this proof lies exactly the matrix version of the variance we have introduced.
As an immediate application of our variance bounds we
obtain stronger versions of B{\"o}ttcher and Wenzel's upper bound.
\end{abstract}

\end{frontmatter}
\section{Variance bounds for a real random variable \label{sec1}}
The variance $\var(X)$ of a random variable $X$
that can assume the real values $x_i$
and does so with probabilities $p_i$
is defined as
\be
\var(X) = \sum_i p_i x_i^2-(\sum_i p_i x_i)^2
= \sum_i p_i(x_i-\sum_j p_j x_j)^2.\label{eq:vardef}
\ee
It is of interest in mathematical statistics to have upper bounds on this variance.
A simple upper bound is given by
\be
\var(X) \le \sum_i x_i^2/2,\label{eq:varb1}
\ee
which follows directly from a much sharper variance bound, due to Murthy and Sethi \cite{ms}.
\begin{lemma}[Murthy-Sethi]
Let $X$ be a real random variable satisfying $m\le X\le M$.
Then $\var(X)\le (M-m)^2/4$.
\end{lemma}
Since $(M-m)^2/4 = (m^2+M^2)/2-(m+M)^2/4 \le (m^2+M^2)/2 \le \sum_i x_i^2/2$,
this bound immediately implies the bound (\ref{eq:varb1}).

\textit{Proof.} The argument, adapted from Muilwijk \cite{muil}, goes as follows.
Some elementary algebra will convince the reader of the following equality:
$$
\var(X) = \sum_i p_i (x_i-m)(x_i-M) + (\mu-m)(M-\mu),
$$
where $\mu=\sum_i p_i x_i$. Because $x_i-m\ge0$ and $x_i-M\le0$, 
the first term is non-positive (while the second is non-negative). Hence
$(\mu-m)(M-\mu)$ is an upper bound on $\var(X)$.
By the arithmetic-geometric inequality,
$$
\sqrt{(\mu-m)(M-\mu)}\le ((\mu-m)+(M-\mu))/2 = (M-m)/2,
$$
and the bound follows.

The inequality is sharp as equality is achieved for a distribution where $X$ is either
$m$ or $M$ with probability $1/2$.
\qed

In this paper we will derive various generalisations of the Murthy-Sethi (MS) bound, and will highlight
its geometric nature. The first generalisation concerns complex-valued random variables
(section \ref{sec:complex}), and this will carry over
in a straightforward way to a matrix generalisation of variance, 
in the special case that the matrix is normal (section \ref{sec:normal}).
Then, in section \ref{sec:nonnormal}, we consider our main objective of a generalisation of variance 
that includes non-normal matrices. 
Along the way we relate these variance bounds to 
the concept of radius of a set of points, and to the new concept we introduce here 
of Cartesian radius of a matrix (not to be confused
with spectral radius, nor with numerical radius).

Before we embark on these generalisations, however, we first describe the seemingly unrelated problem of 
finding sharp bounds on certain norms of a commutator $[X,Y]$ in terms of the norms of $X$ and $Y$
(section \ref{sec:commutator}).
In section \ref{sec:proof}
we give a new proof of a known result and show that at the heart of it lies the concept of variance of a matrix.
The variance bounds we will obtain in this paper can therefore be applied to commutators straight away
and allow us to derive new bounds on norms of commutators.
\section{Notations\label{sec:notations}}
In this paper we are concerned with several kinds of matrix norms.
First of all, as the most general class we'll consider the unitarily invariant (UI) norms, which we denote
using the symbol $|||.|||$. As is well-known, any UI norm of a matrix $X$ 
can be expressed in terms of the singular values of $X$, denoted $\sigma_i(X)$.
As is customary, we assume that singular values are sorted in non-decreasing order. For an $n\times m$ matrix
$X$,
$\sigma_1(X)\ge \sigma_2(X)\ge\ldots\ge\sigma_N(X)\ge0$, with $N=\min(n,m)$.

Special classes of UI norms are the Schatten $p$-norms, the Ky Fan $k$-norms, and the Ky Fan $(p,k)$-norms.
The Schatten $p$-norms are the non-commutative analogues of the $\ell_p$ norms and are defined,
for any $p\ge1$, as
$$
||X||_p := (\trace |X|^p)^{1/p},
$$
where $|X|$ denotes the (left)-modulus of $X$,
$$
|X| := (X^* X)^{1/2}.
$$
In terms of singular values, $||X||_p = (\sum_{i=1}^N \sigma_i(X)^p)^{1/p}$.
For $p=2$, we retrieve the Frobenius norm, also called Hilbert-Schmidt norm,
$$
||X||_2 = \sqrt{\sum_{i=1}^n\sum_{j=1}^m |X_{ij}|^2}.
$$
The Ky Fan $k$-norms are the sums of the $k$ largest singular values,
$$
||X||_{(k)} = \sum_{i=1}^k \sigma_i(X).
$$
Intermediate between these norms are the Ky Fan $(p,k)$-norms \cite{HJII}, which are defined as
$$
||X||_{(k),p} = (\sum_{i=1}^k \sigma_i(X)^p)^{1/p}.
$$

We will use several special matrices repeatedly:
the $n\times n$ identity matrix $\id_n$, or just $\id$ if there is no risk of confusion; 
the standard matrix basis element $e^{ij}$, which has a 1 in position $(i,j)$ and all zeroes elsewhere;
the standard vector basis element $e^i$;
and the Pauli matrices known from quantum physics,
$$
\sigma_x = \twomat{0}{1}{1}{0},
\sigma_y = \twomat{0}{i}{-i}{0},
\sigma_z = \twomat{1}{0}{0}{-1}.
$$
We denote the diagonal matrix with diagonal elements $(x_1,x_2,\ldots,x_n)$ by $\diag(x_1,x_2,\ldots,x_n)$.
Finally, we need the matrix $\diag(1,1,0,\ldots,0)$ so often that we assign it the symbol $F$.

We also use a concept from quantum mechanics called the density matrix.
Disregarding the physical interpretations, we call a matrix a density matrix iff it is
positive semidefinite and has trace 1. This implies that both the vector of eigenvalues and the vector
of diagonal elements (in any orthonormal basis) are formally discrete probability distributions,
being composed of non-negative numbers and summing to 1.
We denote density matrices by lower case greek letters $\rho$ and $\sigma$.
The set of $d\times d$
density matrices is convex and its extremal points are the rank 1 matrices $\psi\psi^*$,
where $\psi$ can be any normalised vector in $\C^d$.
\section{Commutator bounds\label{sec:commutator}}
The commutator of two matrices (or operators) $X$ and $Y$ is defined as $[X,Y]=XY-YX$ and plays
an important role in many branches of mathematics, mathematical physics, quantum physics, and quantum chemistry.
In \cite{bottcher}, B{\"o}ttcher and Wenzel studied the commutator from the following mathematical viewpoint:
fixing the Frobenius norm of $X$ and $Y$, they asked
``How big can the Frobenius norm of the commutator be and how big is it typically?''

By a trivial application of the triangle inequality and H{\"o}lder's inequality one finds
that $||\,[X,Y]\,||_2 \le 2||X||_2 ||Y||_2$.
However, it appears that 2 is not the best constant. It is straightforward to show in the case where $X$ and $Y$
are normal that the best constant is actually $\sqrt{2}$. 
Numerical experiments led B{\"o}ttcher and Wenzel to conjecture that
$\sqrt{2}$ is also the best constant when $X$ and $Y$ are not normal.
Their conjecture can be stated thus:
\begin{theorem}[B{\"o}ttcher and Wenzel]\label{th:frobnorm}
For general complex matrices $X$ and $Y$, and for the Frobenius norm $||.||_2$,
\be
||[X,Y]||_2 \le \sqrt{2} ||X||_2 ||Y||_2.
\ee
The inequality is sharp.
\end{theorem}
We state it here as a theorem because the conjecture has been proved since.

Equality is obtained for $X$ and $Y$ two anti-commuting
Pauli matrices; say $X=\sigma_x$ and $Y=\sigma_z$,
then $[X,Y]=-2i\sigma_y$. This gives
$||[X,Y]||_2=2\sqrt{2}$ and $||X||_2=||Y||_2=\sqrt{2}$.

As already mentioned, the case of normal matrices is rather easy.
For non-normal real $2\times 2$ matrices the proof is also easy, and
Laszlo proved the $3\times 3$ case \cite{laszlo}.
The first proof for the real $n\times n$ case was found by Seak-Weng Vong and Xiao-Qing Jin \cite{VJ} 
and independently by Zhiqin Lu \cite{Lu}. 
Finally, B{\"o}ttcher and Wenzel found a simpler proof \cite{BW} 
that also includes the complex $n\times n$ case.

The empetus behind the present paper was the desire to find an even shorter 
and more conceptual proof, that would
also allow natural generalisations to prove extensions of the theorem.
One can indeed ask for the sharpest constant when $[X,Y]$, $X$ and $Y$ are compared
in terms of different Schatten norms. That is:
\begin{problem}
Let $X$ and $Y$ be general square matrices, and $p,q,r\ge1$,
such that $1/p\le 1/q+1/r$ holds.
What is the smallest value of $c$ such that
$$
||[X,Y]||_p \le c ||X||_q ||Y||_r
$$
holds?
We will denote this smallest $c$ by $c_{p,q,r}$.
\end{problem}
The restriction $1/p\le 1/q+1/r$ is necessary. When it is not satisfied, $c$ is dimension
dependent, just as in the case of H\"older's inequality.
To see this, take two fixed non-zero $X$ and $Y$ for which $||[X,Y]||_p \le c ||X||_q ||Y||_r$ holds, with
some predetermined finite value of $c$.
Then replace $X$ and $Y$ by $X\otimes\id_D$ and $Y\otimes \id_D$, with large value of $D$.
The left-hand side of the inequality is thus multiplied by $D^{1/p}$, 
while the right-hand side is multiplied by $D^{1/q+1/r}$.
If $1/p\ge 1/q+1/r$, the left-hand side grows faster with $D$ than the right-hand side, and for some
large enough $D$, the inequality will be violated for any initial choice of $c$. 
Thus, if $1/p\ge 1/q+1/r$, there is no finite $c$ for which the inequality holds universally,
and henceforth we only consider the case when $1/p\le 1/q+1/r$.

Numerical experiments have led us to conjecture:
\begin{conjecture} 
\label{conj:gen}
For the restricted case $p=q$,
\be
c_{p,q,r}=c_{p,p,r}=2^{\max(1/p,1-1/p,1-1/r)}.
\ee
\end{conjecture}
By taking special examples of $X$ and $Y$ we can calculate lower bounds on the constant $c_{p,q,r}$.
This allows us to check that the conjectured bounds would be sharp. 
\begin{enumerate}
\item The two anti-commuting Pauli matrices $X=\sigma_x$ and $Y=\sigma_z$ give 
$||[X,Y]||_p=2^{1+1/p}$ and $||X||_q=2^{1/q}$ and $||Y||_r=2^{1/r}$, hence
$$
c\ge 2^{1+1/p-1/q-1/r}.
$$
\item The choices $X=e^{12}$ and $Y=e^{21}$ give $[X,Y]=\sigma_z$, hence $||\,[X,Y]\,||_p=2^{1/p}$ and
$||X||_q=||Y||_r=1$, thus 
$$
c\ge 2^{1/p}.
$$
\item The two anti-commuting matrices
$$
X=\twomat{\sqrt{2}}{-2-\sqrt{2}}{2-\sqrt{2}}{-\sqrt{2}}/4,\quad Y=\twomat{1}{1}{1}{-1}/\sqrt{2}.
$$
Since $X$ is rank 1 and $Y$ is unitary,
both $XY$ and $[X,Y]$ are rank 1.
The singular values of $X$ are $(1,0)$,
those of $Y$ are $(1,1)$, and those of $[X,Y]$ are $(2,0)$.
Thus $||\,[X,Y]\,||_p=2$, $||X||_q=1$ and $||Y||_r=2^{1/r}$. Of course, $X$ and $Y$ can be swapped.
This gives 
$$
c\ge 2^{1-1/r}, \qquad
c\ge 2^{1-1/q}.
$$
\end{enumerate}

Note that in \cite{BW} the special case $p=q=r$ of this conjecture has already appeared
(eq.\ (28) in \cite{BW}), which subsequently has been proven by Wenzel, using complex interpolation
(Riesz-Thorin) methods \cite{wentzel}.
The methods investigated in the present paper will allow us to
establish the conjecture for $p=q=2$ and general $r$.

Furthermore, the special case $p=q=r=\infty$, has been proven a long time ago by Stampfli \cite{stampfli}
in the broader setting of operator algebras. The goal there was to study the operator norm of
the operator $D_Y: X\mapsto [X,Y]$ in terms of the operator norm of $Y$. 
In that sense, our conjecture relates the norm of $D_Y$ acting on Schatten class $L_p$ to the Schatten $r$ norm
of $Y$.

In the regime where $p$, $q$ and $r$ satisfy $1/p=1/q+1/r$, the best constant is trivial and equal to 2.
In fact, we have the more general theorem for all UI norms:
\begin{theorem}\label{th:qnorm}
For general complex matrices $X$ and $Y$, and for any UI norm $|||.|||$,
\be
|||\,[X,Y]\,||| \le 2 ||| \,|X|^s\, |||^{1/s}\,\, ||| \,|Y|^t\, |||^{1/t}.
\ee
\end{theorem}
\textit{Proof.}
This just follows from a combination of 
the triangle inequality and the fact that $|||XY|||=|||YX|||$ for any UI norm, 
with H\"older's inequality for UI norms:
$$
|||\,[X,Y]\,||| = |||XY-YX||| \le |||XY|||+|||YX||| = 2|||XY|||,
$$
and
$$
|||XY||| \le ||| \,|X|^s\, |||^{1/s}\,\, ||| \,|Y|^t\, |||^{1/t},
$$
for $s$ and $t$ satisfying $s>1$ and $1/s+1/t=1$
(\cite{bhatia}, Corollary IV.2.6).
\qed

In spite of the triviality of the proof, the factor 2 is the best constant. 
Indeed, equality is obtained for $X$ and $Y$ two anti-commuting
Pauli matrices.

Applied to Schatten $p$-norms, this gives the special case mentioned above
\be
||[X,Y]||_p \le 2 ||X||_{q} ||Y||_{r},
\ee
for $1/p=1/q+1/r$.

In the following section we present our proof of Theorem \ref{th:frobnorm}, and a certain expression
obtained halfway through it will allow us to make contact with our main object of interest, namely
the variance bounds mentioned at the beginning.
\section{A New, Shorter Proof of Theorem \ref{th:frobnorm}\label{sec:proof}}
One easily checks the following:
\beas
||XY-YX||_2^2     &=& \trace[XYY^*X^*-XYX^*Y^*-YXY^*X^*+YXX^*Y^*] \\
                  &=& \trace[X^*XYY^*-XYX^*Y^*-YXY^*X^*+XX^*Y^*Y] \\
||X^*Y+YX^*||_2^2 &=& \trace[YX^*Y^*X+YX^*XY^*+X^*YXY^*+X^*YY^*X] \\
                  &=& \trace[X^*XY^*Y+XYX^*Y^*+YXY^*X^*+XX^*YY^*].
\eeas
Taking the sum yields
\bea
\lefteqn{||XY-YX||_2^2 + ||X^*Y+YX^*||_2^2} \nonumber\\
&=& \trace[X^*XYY^*+XX^*Y^*Y+X^*XY^*Y+XX^*YY^*] \nonumber\\
&=& \trace(X^*X+XX^*)(Y^*Y+YY^*).\label{eq:prf1}
\eea
By the Cauchy-Schwarz inequality,
\bea
|\trace[Y(X^*X+XX^*)]| &=&
|\trace[(YX^*+X^*Y)X]| \nonumber\\
&\le& ||YX^*+X^*Y||_2\,\,||X||_2.\label{eq:prf2}
\eea
Combining (\ref{eq:prf1}) and (\ref{eq:prf2}) then gives
\beas
||XY-YX||_2^2 &\le& \trace(X^*X+XX^*)(Y^*Y+YY^*) \\
&& -|\trace[Y(X^*X+XX^*)]|^2/||X||_2^2.
\eeas
Introducing the matrix $\rho=(X^*X+XX^*)/(2||X||_2^2)$, this can be expressed as
\be
||XY-YX||_2^2 \le 4||X||_2^2 \left(\trace[\rho(Y^*Y+YY^*)/2]-|\trace[\rho Y]|^2\right).
\ee
Note that $\rho$ is positive semi-definite and has trace 1 and is formally a density matrix.
The quantity $\trace[\rho(Y^*Y+YY^*)/2]-|\trace[\rho Y]|^2$ appearing here is reminiscent
of the variance of a random variable, with $\rho$ taking over the role of a probability distribution.

To make the connection even more obvious,
consider now the Cartesian decomposition $Y=A+iB$, where $A$ and $B$ are Hermitian.
One checks that $(Y^*Y+YY^*)/2=A^2+B^2$.
Therefore,
$$
\trace[\rho(Y^*Y+YY^*)/2]-|\trace[\rho Y]|^2 = \trace \rho(A^2+B^2)-(\trace \rho A)^2-(\trace \rho B)^2,
$$
which is a sum of terms in $A$ and in $B$ separately.
We now need to show that the right-hand side
is bounded above by $||Y||_2^2/2 = (||A||_2^2+||B||_2^2)/2$.

This would follow if
$\trace \rho A^2-(\trace \rho A)^2 \le ||A||_2^2/2$
for all Hermitian $A$.
We can prove this by passing to a basis in which $A$ is diagonal, so let's put
$A=\diag(a_1,\ldots,a_d)$ and let us denote the diagonal
elements of $\rho$ in that basis by $p_i$. As the $p_i$ are non-negative and add up to 1,
they form a probability distribution. The quantity $\trace \rho A^2-|\trace \rho A|^2$ then becomes
$$
\sum_i p_i a_i^2-(\sum_i p_i a_i)^2.
$$
This is the variance $\var(A)$ of a random variable $A$
that can assume the values $a_i$
and does so with probabilities $p_i$.
Applying the variance bound
$$
\var(A) \le \sum_i a_i^2/2=||A||_2^2/2
$$
then proves the required statement.
\qed

In the remainder of the paper we study the quantity
$$\trace[\rho(Y^*Y+YY^*)/2]-|\trace[\rho Y]|^2,
$$
which
can be seen as a generalisation of the variance of a random variable to the matrix (quantum) case.
We derive sharp bounds on this generalised variance, which directly lead to sharper bounds on the 2-norm
of a commutator, and which allow us to prove a special case of
our conjecture about commutator norms.
\section{Variance of a complex random variable\label{sec:complex}}
First of all, we formally define the variance of a complex random variable.
To do so, we replace squares in (\ref{eq:vardef}) by modulus square, whether this makes
statistical sense or not.
For $x_i\in \C$:
\be
\var(X) = \sum_i p_i |x_i|^2-|\sum_i p_i x_i|^2
= \sum_i p_i|x_i-\sum_j p_j x_j|^2.\label{eq:vardefC}
\ee
In statistical terms, this corresponds to the trace of the covariance matrix when considering real and imaginary
part of $X$ as two random variables.

Our first result, proven below, is a straight generalisation of the MS bound to the complex case.
\begin{theorem}\label{th:radiusC}
For a random variable $X$ assuming complex values $x_i$, the largest possible variance obeys
\be
\max_{\bm{p}} \sum_i p_i|x_i-\sum_j p_j x_j|^2
= \min_{y\in\C} \max_i |x_i-y|^2.
\ee
\end{theorem}
The right-hand side can be interpreted in the context of Euclidean planar geometry applied to the
complex plane (with the modulus acting as Euclidean norm).
\begin{definition}
The \emph{radius} of a set of points $\cX=\{x_i\}$ in the Euclidean plane, denoted $r(\cX)$, is
the radius of the smallest circle circumscribing $\cX$.
The \emph{center} of $\cX$ is the center of that circle.
\end{definition}
Theorem \ref{th:radiusC} thus says that the variance of $X$ taking values in $\cX$
is bounded above by the square of $r(\cX)$.

Some obvious properties of the radius of a set
are that it is invariant under global translations, rotations and reflections.
It is also homogeneous of degree 1.
\subsection{Proof of Theorem \ref{th:radiusC}}
A simple observation will be important for the proof.
Let $C$ be the smallest circumscribing circle of $\cX$ and let $c$ be its center.
Let $\cX'$ be the subset of points of $\cX$ that lie on $C$. There must at least be two such points,
for if it contained only 1 point
a smaller circle could be found by moving the center towards that point.
Then, $\forall x\in \cX', |x-c|=r(\cX)$, and for all other $x$, $|x-c|<r(\cX)$ strictly.
This means that for a new point $c'$ close enough to $c$, $\max_i |x_i-c'|$ is obtained for $x_i$ on $C$.

\begin{figure}[ht]
\begin{center}
\includegraphics[width=9cm]{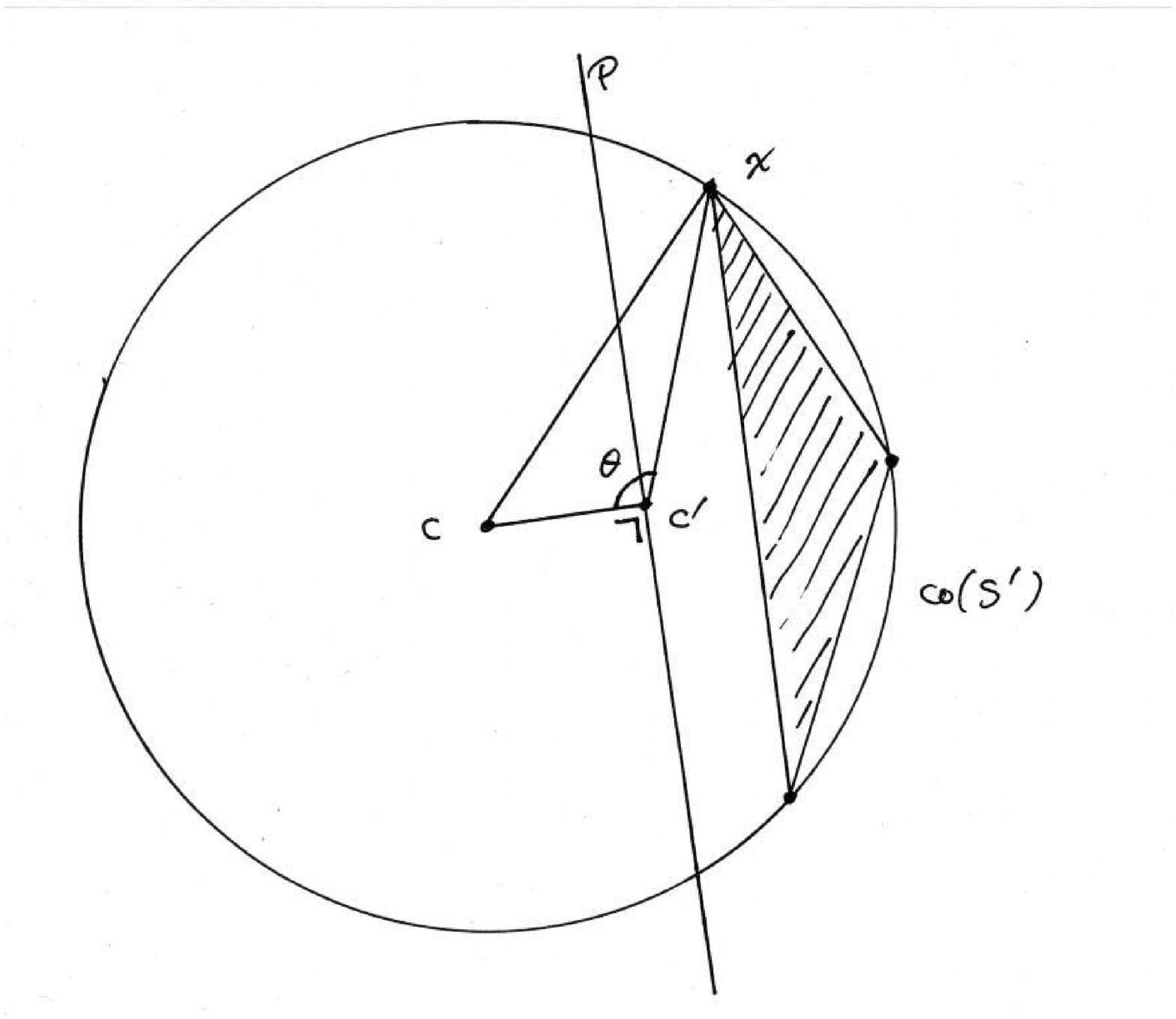}
\caption{
\label{fig:center}
}
\end{center}
\end{figure}

\begin{lemma}\label{lem:center}
The center of a set $\cX$ endowed with a Euclidean metric is contained in the convex hull
of the points of $\cX'$.
\end{lemma}
\textit{Proof.}
Suppose, to the contrary, that $c$ lies outside the convex hull of $\cX'$.
By Minkowski's separating hyperplane theorem there must then be a hyperplane
$P$ such that $c$ is strictly on one side, while
all points of $\cX'$ are strictly on the other side of $P$.
Let $c'$ be the orthogonal projection of $c$ on $P$ and let $x$ be any point in $\cX'$.
From the geometry follows that the triangle $c,c',x$
has an obtuse angle at $c'$ (see Figure \ref{fig:center}).
By the cosine rule one then sees that every point $x\in \cX'$ is strictly closer
to any point on the open line segment $]cc'[$ than to $c$, violating the assumption
that $c$ is the center of $\cX$.
\qed

\textit{Proof of Theorem \ref{th:radiusC}.}
We start with the expression
$ \sum_i p_i|x_i-\sum_j q_j x_j|^2$, where $\bm{p}$ and $\bm{q}$ are two probability distributions.
As any average of real quantities is bounded above by the maximum, we can replace the outer average and get
$$
\sum_i p_i|x_i-\sum_j q_j x_j|^2
\le
\max_i |x_i-\sum_j q_j x_j|^2.
$$
This is true for any $\bm{q}$. Hence, the inequality remains if both sides are minimised over $\bm{q}$.

The minimisation of the right-hand side,
$\min_{\bm{q}} \max_i |x_i-\sum_j q_j x_j|^2$, is almost the right-hand side of the Theorem,
but with the minimisation
over any complex value $y$ replaced by a minimisation over the convex hull of the set of $x_i$.
From lemma \ref{lem:center}, however, we see that the optimal $y$ will be
within that convex hull, so that both minimisations
must yield the same value.

We will now show that the left-hand side is minimal for $\bm{q}$ equal to $\bm{p}$, in which case the value
is equal to the variance.
Let $\mu_P=\sum_i p_i x_i$ and $\mu_Q=\sum_i q_i x_i$.
Obviously, $|\mu_Q-\mu_P|^2\ge 0$. Thus
$|\mu_Q|^2 - 2\Re \overline{\mu_Q}\mu_P \ge |\mu_P|^2 - 2\Re \overline{\mu_P}\mu_P$.
From this it follows immediately that
$\sum_i p_i|x_i-\mu_Q|^2 \ge \sum_i p_i |x_i-\mu_P|^2$, for all $\bm{q}$.

To show that equality holds, take $p_j$ such that $\sum_{j=1}^m p_j x_j = y^*$, where only points $x_j$
in $\cX'$ contribute. This is possible because, by lemma \ref{lem:center}, $y^*$
is in the convex hull of those points.
\qed
\subsection{Relation between radius and vector norms}
The radius is not a norm, because it is not convex. Nevertheless,
our next two results draw the connection between the radius and permutation invariant (PI) vector norms.
First we show that the radius is bounded above by one half the value of a specific PI vector norm and then we
derive from that how it relates to all other PI vector norms, giving best constants for each.
We introduce some notation, borrowed from the theory of majorisation:
let $|\cX|^{\downarrow,k}$ be the $k$-th largest value among the moduli of $\cX=\{x_i\}_i$.
We freely consider $\cX$ either as a set of $d$
points in $\C$ or as a vector in $\C^d$.
We also use the shorthand $\cX-z = \{x_i-z\}_i$.

The central statement is that the maximum in the definition of $r(\cX)$ can be replaced by means of the
largest and the second largest value.
\begin{theorem}\label{th:rkf2}
For any set of complex values $\cX=\{x_i\}_{i=1}^d$, and any $p\ge 1$,
\be
r(\cX) = \min_{z\in \C} \left(\frac{
(|\cX-z|^{\downarrow,1})^p+
(|\cX-z|^{\downarrow,2})^p}{2}\right)^{1/p}.\label{eq:rkf2}
\ee
\end{theorem}
By putting $z=0$, we then immediately get:
\begin{corollary}
For any set of complex values $\cX=\{x_i\}_{i=1}^d$, and any $p>1$,
\be
r(\cX) \le \left(((|\cX|^{\downarrow,1})^p+(|\cX|^{\downarrow,2})^p)/2\right)^{1/p}.
\ee
\end{corollary}
The relation to all other PI norms then follows from:
\begin{theorem}\label{th:kyfanvec}
For all permutation invariant norms $||.||$ on $\C^d$,
\be
(|\cX|^{\downarrow,1}+|\cX|^{\downarrow,2})/2
\le \frac{||\cX||}{||F||}
\le \max\left(||\cX||_\infty, ||\cX||_1/2 \right).
\ee
\end{theorem}
Because the last theorem is easily generalised to matrices in terms of the Ky Fan $k$-norm $||X||_{(k)}$,
we will prove it for matrices straight away.
\begin{theorem}\label{th:kyfan}
For all unitarily invariant norms $|||.|||$ on $\M(\C^d)$,
\be
||X||_{(2)}/2 \le \frac{|||X|||}{|||F|||} \le \max\left(||X||_\infty, ||X||_1/2 \right).
\ee
\end{theorem}
Theorem \ref{th:kyfanvec} follows by setting $X=\diag(x_i)$.

\textit{Proof.} 
We start with the lower bound.
Note first that
$$
||X||_{(2)}/2 = ||X||_{(2)}/||F||_{(2)}.
$$
We wish to prove that of all UI norms, the $2$nd Ky Fan norm minimises the ratio $||X||/||F||$.

Every unitarily invariant norm $|||.|||$ can be defined as
(\cite{HJII}, Theorem 3.5.5)
$$
|||X||| = \max\{ \sum_i \alpha^{\downarrow}_i \sigma_i(X): \alpha\in N_{|||.|||}\},
$$
where $N_{|||.|||}$ is a compact subset of $\R_+^d$ specific to that norm and
$\sigma_i(X)$ are the singular values of $X$.
Minimising over all UI norms thus amounts to minimising over all compact sets $N$.
In particular, minimising the ratio $|||X|||/|||F|||$ amounts to minimising over all compact sets $N$
whose associated norm obeys the constraint $|||F|||=1$, that is
$\max_{\alpha\in N} \alpha^{\downarrow}_1+\alpha^{\downarrow}_2=1$.
Thus
\beas
\min_{|||.|||} |||X|||/|||F|||
&=& \min_N \max_{\alpha\in N}\{ \sum_i \alpha^{\downarrow}_i \sigma_i(X):
\alpha^{\downarrow}_1+\alpha^{\downarrow}_2\le 1\} \\
&=& \min_N \max_{\alpha\in N}\{ \sum_i \alpha^{\downarrow}_i \sigma_i(X):
\alpha^{\downarrow}_1+\alpha^{\downarrow}_2 = 1\} \\
&=& \min_\alpha \{\alpha_1^{\downarrow}\sigma_1(X) + \alpha_2^{\downarrow}\sigma_2(X):
\alpha^{\downarrow}_1+\alpha^{\downarrow}_2 = 1\} \\
&=& (\sigma_1(X)+\sigma_2(X))/2,
\eeas
which proves the lower bound.

For the upper bound, we similarly have
$$
\max_{|||.|||} |||X|||/|||F|||
= \max_N \max_{\alpha\in N}\{ \sum_i \alpha^{\downarrow}_i \sigma_i(X):
\alpha^{\downarrow}_1+\alpha^{\downarrow}_2\le 1\}.
$$
Thus elements $\alpha^{\downarrow}_3$ and beyond must be as large as possible, which means they should
be equal to $\alpha^{\downarrow}_2$. The maximisation then reduces to a maximisation over $\alpha^\downarrow_1=:a$,
where $1/2\le a\le 1$,
\beas
\max_{|||.|||} |||X|||/|||F|||
&=& \max_{1/2\le a\le 1} a\sigma_1(X)+(1-a)\sum_{k=2}^d\sigma_k(X) \\
&=& \max_{1/2\le a\le 1} (2a-1)\sigma_1(X) + (1-a)\sum_{k=1}^d\sigma_k(X) \\
&=& \max_{1/2\le a\le 1} (2a-1)||X||_\infty + (1-a)||X||_1.
\eeas
The maximum is attained in one of the extreme points, $a=1/2$ or $a=1$, hence
$$
\max_{|||.|||} |||X|||/|||F|||
=\max(||X||_1 /2,||X||_\infty).
$$
\qed

\begin{corollary}\label{cor:kyfan}
For all $p\ge p_0$,
\be
2^{-1/p_0} ||\,\,|X|^{p_0}\,\,||_{(2)}^{1/p_0}
\le \frac{||X||_p}{||F||_p}
\le \max\left(2^{-1/p_0} ||X||_{p_0} , ||X||_{\infty}\right).
\ee
\end{corollary}
Note that the norm in the left hand side is the Ky Fan $(p,k)$-norm $||X||_{(2),p_0}$.

\textit{Proof.}
Apply theorem \ref{th:kyfan} to $|X|^{p_0}$ and note that $||\,\,|X|^{p_0}\,\,||_q = ||X||_{p_0q}^{p_0}$.
\qed

\subsection{Proof of Theorem \ref{th:rkf2}}
To prove Theorem \ref{th:rkf2}, we need a lemma.
\begin{lemma}\label{lem:mid}
Consider a polygon $P$. Let $P'$ be the polygon whose vertices are the midpoints of the edges of $P$.
Then the center of the smallest circle that circumscribes $P$ is in $P'$.
\end{lemma}
\textit{Proof.}
Consider first the simplest case that $P$ is a triangle $ABC$. By a well-known and easily proven
geometrical theorem,
the center of a circumscribing sphere containing all three points of the triangle
is equal to the intersection $D$ of the bisectors of the triangle's edges.
By definition, these bisectors pass through the midpoints of the edges of $P$, which are the vertices of $P'$.
By inspection one sees that if $D$ lies in $P=ABC$, it must also lie in $P'$.
The two possible cases are illustrated in Figure 1.
\begin{figure}[ht]
\includegraphics[width=13cm]{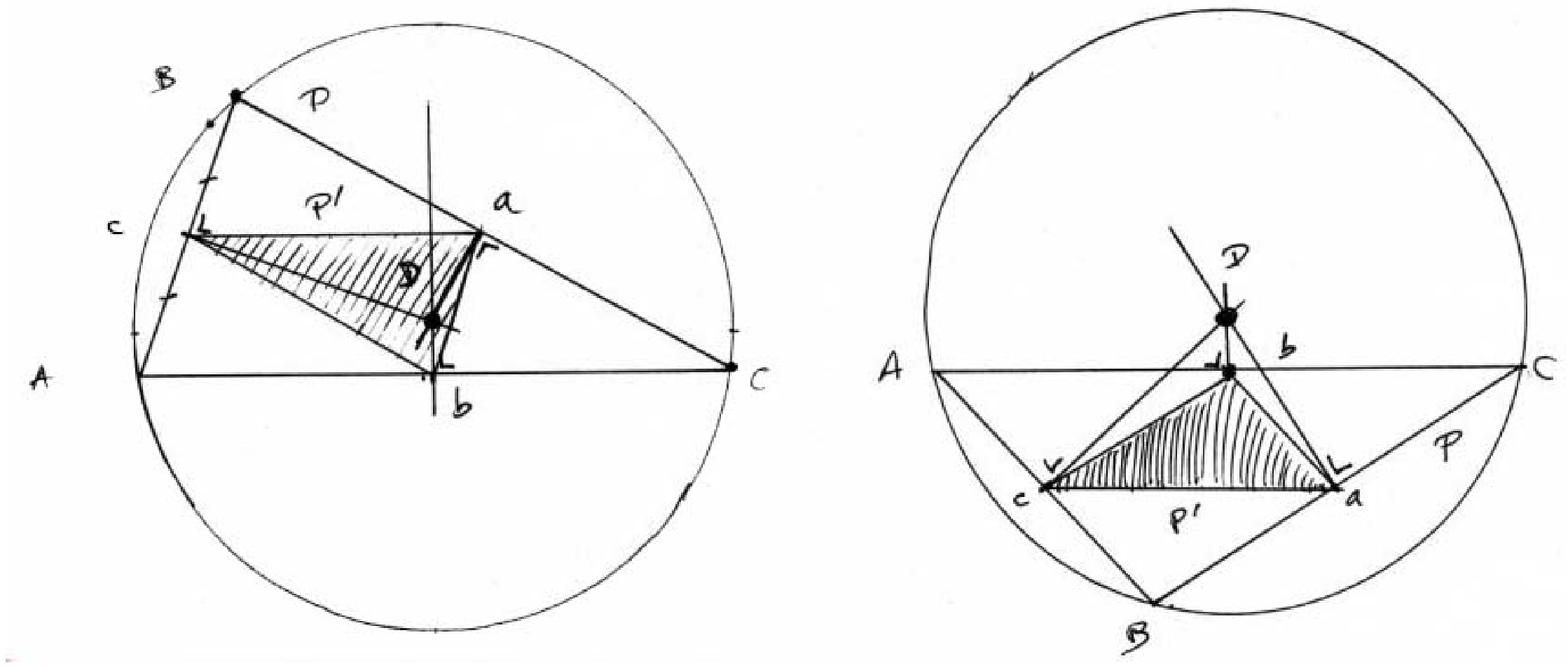}
\caption{
\label{fig:rkf2}
}
\end{figure}

If $P$ is a general polygon, it can be subdivided into one or more non-overlapping triangles.
According to Lemma \ref{lem:center}, the center $D$
of the smallest circle circumscribing $P$ is in $P$. Therefore, $D$ is in
one of those triangles; call it $ABC$. By the above argument, $D$ is also in the triangle of midpoints
of $ABC$. This triangle of midpoints is a subset of the polygon $P'$ of midpoints. Hence $D$ is in $P'$ too.
\qed

\textit{Proof of Theorem \ref{th:rkf2}.}
We will show that the optimal $z$
in the RHS of (\ref{eq:rkf2}) is equal to $y^*$, the optimal $y$ in $r(\cX)=\min_y \max_i |x_i-y|$.

We relabel the points of $\cX$
so that $x_1,x_2,\ldots,x_m$ are the points in $\cX'$, i.e.\ they are the points on the
smallest circumscribing circle around $\cX$,
which has center in $y^*$. There must be at least two points in $\cX'$.
Thus, $|\cX-y^*|^{\downarrow,1}=|\cX-y^*|^{\downarrow,2}$, so that the RHS of
(\ref{eq:rkf2}) is equal to its LHS in the point $z=y^*$.

We must show that the RHS is minimal in $z=y^*$.
We begin by pointing out that the RHS is a norm of $\cX-z$ and
hence a convex function of $z$ (see, e.g.\, \cite{bhatia}, Example IV.1.4).
Thus, this function has a single local minimum, which automatically is the global minimum.
To find out whether $z=y^*$ is indeed the global minimum, it suffices to check whether it is a local minimum.
We'll do so by perturbing $z$ by an infinitesimal amount: $z=y^*+t\Delta$.

If $t$ is small enough, the only contributions to the derivative of
$(|\cX-z|^{\downarrow,1})^p+(|\cX-z|^{\downarrow,2})^p$
come from the derivatives of $|x_1-z|$, $|x_2-z|$, \ldots, $|x_m-z|$.
More precisely, only the two largest of these derivatives contribute.
The derivative of $|x_i-z|^p$ w.r.t.\ $t$ in $t=0$ is
$p|x_i-y^*|^{p-1}$ (a constant factor for $x_i\in \cX'$) times the
derivative of $|x_i-z|$ w.r.t.\ $t$ in $t=0$, which is $-\langle\Delta,x_i-y^*\rangle$.

To show that $y^*$ is a local minimum, we have to show
that for any $\Delta$ the sum of the two largest derivatives is non-negative.
This means that, for any $\Delta$, there exist distinct $i$ and $j$ such that
$-\langle\Delta,x_i-y^*\rangle-\langle\Delta,x_j-y^*\rangle\ge0$, i.e.\
$\langle(-\Delta),(x_i+x_j)/2-y^*\rangle\ge0$. Now, the set of points $(x_i+x_j)/2$ contains the midpoints
of the edges of the polygon $P$ with vertices $x_1,x_2,\ldots,x_m$.
By Lemma \ref{lem:mid} the polygon $P'$ whose vertices are these midpoints contains the center $y^*$
of the circle. Therefore, for any direction $\Delta$ there will be some midpoint $(x_i+x_j)/2$
such that $\langle(-\Delta),(x_i+x_j)/2-y^*\rangle\ge0$. This shows that, indeed, $z=y^*$ is a local minimum.
\qed

This proof relies heavily on planar geometry. It would be interesting
to find an entirely algebraic proof.

\subsection{Main Result}
Combining all results obtained so far yields the main theorem of this section:
\begin{theorem}
For a complex valued random variable $X$, taking values in the discrete set $\cX=\{x_i\}_{i=1}^d$,
and for any PI vector norm $||.||$,
\be
\sqrt{\var(X)} \le r(\cX) \le ||\cX||_{(2)}/2 \le ||\cX||/||F||,
\ee
where $\cX$ has been interpreted as a vector in $\C^d$, and $F=(1,1,0,\ldots,0)$.
\end{theorem}

\textbf{Remark.}
Many other generalisations are possible of the concepts introduced here.
In the above we've considered complex valued $X$,
with norm given by the complex modulus. This is isomorphic to vectors in $\R^2$, endowed with the
Euclidean 2-norm. We can more generally consider $X$ whose values are in $\ell_p$,
or even in the Schatten class $L_p$.
\section{Quantum Variance of Normal Matrices\label{sec:normal}}
In this section we consider variance bounds in the matrix setting, where
probability distributions are replaced by density matrices.
This leads to the following
definition of the variance of a normal matrix $X$:
\begin{definition}
The \emph{quantum variance} of a \emph{normal} matrix $X\in \M_d(\C)$ w.r.t.\ the 
density matrix $\rho$ is given by
\end{definition}
\be
\var(X) = \trace[\rho |X|^2] - |\trace[\rho X]|^2
=\trace[\rho |X-\trace[\rho X]\id_d|^2].
\ee
Here, $|.|$ stands for the matrix modulus defined by $|X|=(X^*X)^{1/2}$,
and $\id_d$ is the $d\times d$ identity matrix.

\textbf{Remark.} In the mathematical physics literature a more general version of the variance can be found,
based on unital completely positive
maps $\Phi$ \cite{bhatiadavis}, where the variance is operator-valued.
Our definition here corresponds to the choice $\Phi(X)=\trace[\rho X]$,
yielding a scalar-valued variance.

This definition is a straightforward generalisation of the classical variance for complex scalar variables.
Since $X$ is normal, it can be diagonalised by a unitary conjugation.
Inserting $X=U\Lambda U^*$ in the definition,
with $\Lambda=\diag(\lambda_1,\ldots,\lambda_d)$ complex, yields
$$
\var(X) = \sum_i p_i |\lambda_i|^2-|\sum_j p_j\lambda_j|^2,
$$
where $p_i$ is the diagonal element $(U^*\rho U)_{ii}$. Therefore, if $\rho$ can be any density matrix,
$\bm{p}=(p_1,\ldots,p_d)$ can be any probability distribution.
It follows that the variance bounds obtained for complex variables carry over wholesale to normal matrices,
by applying them to the spectrum of the normal matrix.

In particular, the radius $r(X)$ of a normal matrix is the radius of its spectrum. Note, however, that the term
spectral radius is already in use and denotes the radius of the smallest circumscribing
circle with center at the origin. One can easily show that the spectral radius of a normal matrix is an upper
bound on the radius of its spectrum.

The main result of the last section becomes:
\begin{theorem}\label{th:normalmain}
For a normal $d\times d$ matrix $X$
and for any UI vector norm $|||.|||$,
\be
\sqrt{\var(X)} \le r(X) \le ||X||_{(2)}/2 \le |||X|||/|||F|||,\mbox{ e.g. } 2^{-1/p}||X||_p.
\ee
\end{theorem}
Using this theorem, B{\"o}ttcher and Wenzel's theorem can already 
be strenghtened in the specific case of normal $Y$, by combining the statement obtained
halfway through its proof with theorem \ref{th:normalmain}.
For normal $Y$, and all $p\ge 1$,
\be
||\,\,[X,Y]\,\,||_2 \le 2 ||X||_2 \,\,r(Y) \le ||X||_2 ||Y||_{(2)} 
\le 2^{1-1/p} ||X||_2 ||Y||_p.
\ee
\section{Quantum Variance of Non-Normal Matrices\label{sec:nonnormal}}
We will now investigate the general case, of quantum variance of a non-normal matrix.
In this case the left modulus and right modulus of $X$,
$(X^*X)^{1/2}$ and $(XX^*)^{1/2}$, are no longer the same.
Therefore, there are many possible distinct extensions of the expression $\trace[\rho |X|^2]$.
One is $\trace[\rho X^* X]$, another is $\trace[\rho XX^*]$,
and we'll also consider the mean of the two, $\trace[\rho(X^*X+XX^*)]/2$, 
which featured prominently in our proof
of Theorem \ref{th:frobnorm}.

For that reason we need a name for the expression $((X^*X+XX^*)/2)^{1/2}$, and we have chosen to call
it the \textit{Cartesian modulus}. One observes that
in terms of the Cartesian decomposition of $X$, $X=A+iB$ with $A$ and $B$ Hermitian, the Cartesian modulus
reduces to the pleasing form $(A^2+B^2)^{1/2}$.

For convenience, we'll denote the three corresponding moduli by $|.|_*$, each with a different subscript:
\bea
|X|_L &:=& (X^* X)^{1/2},\\
|X|_R &:=& (X X^*)^{1/2},\\
|X|_C &:=& \sqrt{(|X|_L^2+|X|_R^2)/2}.
\eea
Note that $|||\,|X|_L\,||| = |||X|||$ for any UI norm and
the same holds for the right modulus. For the Cartesian modulus this is no longer true, but we do have
the following inequalities for Schatten $p$-norms
obtained by Bhatia and Kittaneh (\cite{zhanbook}, eqns (3.38) and (3.39)):
for $p\ge 2$,
\be
|| \,|X|_C\,||_p \le ||X||_p  \le 2^{1/2-1/p} ||\,|X|_C\,||_p, \label{eq:bk}
\ee
while the reversed inequalities hold for $1\le p\le 2$.
More fundamental is the following inequality
for the Ky Fan $(p,k)$-norms with $p=2$
\be
|| \,|X|_C\,||_{(k),2} \le ||X||_{(k),2}, \label{eq:cartkf2}
\ee
for any $k$
(which in \cite{zhanbook} is phrased as a majorisation statement; see its eq.\ (3.31)).

Each modulus builds a different variance, which we'll distinguish by the corresponding subscript too.
Thus
\be
\var_*(X) = \trace[\rho |X-\trace[\rho X]\id_d|_*^2]
\ee
where $*$ stands for $L$, $R$ or $C$.
It is easily checked that
in each case, the variance satisfies the relation
\be
\var_*(X) = \trace[\rho |X|_*^2] - |\trace[\rho X]|^2.
\ee

We next show how to generalise theorem \ref{th:radiusC} to the non-normal matrix case.
In the proof we need the \textit{numerical range} $W(X)$ of a matrix $X$ \cite{HJII}:
$W(X) = \{\psi X \psi^*: \psi\in\C^d, ||\psi||=1\}$.
By the Toeplitz-Hausdorff theorem, $W(X)$ is a convex set. It can therefore be redefined in terms
of density matrices as
\be
W(X) = \{\trace[\rho X]: \rho\ge0, \trace\rho=1\}.
\ee
Henceforth, 
we use the shorthand $\max_\rho$ or $\min_\rho$ to denote maximisation and minimisation over all
possible density matrices $\rho$.
\begin{theorem}\label{th:radiusM}
For a non-normal $n\times n$ matrix $X$,
\be
\var_*(X) \le \max_{\rho} \trace[\rho|X-\trace[\rho X]\id_d|_*^2]
= \min_{y\in\C} || \,\, |X-y\id|_*^2\,\, ||_{\infty}.\label{eq:radiusM}
\ee
Furthermore, the maximisation over $\rho$ can be restricted to density matrices of rank 1, of the form $\psi\psi^*$,
with $\psi$ a normalised vector in $\C^N$.
\end{theorem}
\textit{Proof.}
The proof proceeds in a similar way as in the complex variable case.
The bivariate function
$$
(\rho,\sigma)\mapsto f(\rho,\sigma)=\trace[\rho |X-\trace[\sigma X]\id|_*^2]
$$
satisfies the following properties:
its domains are compact convex sets (being the set of all density matrices),
the function is convex in $\sigma$ for all $\rho$,
concave (linear, in fact) in $\rho$ for all $\sigma$,
and continuous in both $\rho$ and $\sigma$.
All conditions of Kakutani's minimax theorem \cite{kak} are therefore fulfilled, hence
in the minimax expression $\min_\sigma \max_\rho f(\rho,\sigma)$
the minimisation over $\sigma$ and maximisation over $\rho$ can be freely interchanged.

One easily verifies that
\beas
\lefteqn{
\trace[\rho |X-\trace[\sigma X]\id|_R^2]-
\trace[\rho |X-\trace[\rho X]\id|_R^2]} \\
&=& |\trace[\sigma X]-\trace[\rho X]|^2 \\
&\ge& 0,
\eeas
so that
the minimum of $\trace[\rho |X-\trace[\sigma X]\id|_R^2]$
over $\sigma$ is obtained for $\sigma=\rho$.
The same is obviously true for the left modulus, and it also holds
for the Cartesian modulus since $|.|_C^2 = (|.|_L^2+|.|_R^2)/2$.

Therefore, we get the following chain of equalities:
\beas
\max_\rho \trace[\rho |X-\trace[\rho X]\id|_*^2]
&=& \max_\rho \min_\sigma \trace[\rho |X-\trace[\sigma X]\id|_*^2] \\
&=& \min_\sigma \max_\rho \trace[\rho |X-\trace[\sigma X]\id|_*^2] \qquad (*)\\
&=& \min_\sigma ||\, |X-\trace[\sigma X]\id|_*^2\,||_\infty \\
&=& \min_{y\in W(X)} ||\, |X-y\id|_*^2\,||_\infty \\
&=& \min_{y\in\C} ||\, |X-y\id|_*^2\,||_\infty.
\eeas
In the third line we used the Rayleigh-Ritz characterisation of the largest eigenvalue of a Hermitian matrix.
In the last line we could remove the constraint $y\in W(X)$ because of the fact,
proven in lemma \ref{lem:centerW} below, that the optimal $y$
in $\min_{y\in\C}|| \,\,|X-y\id|_*^2\,\,||_\infty$ is automatically in $W(X)$.

To prove the final statement of the theorem, we note that in (*) the maximisation over $\rho$
can be restricted to $\rho$ that have rank 1.
Furthermore, the minimisation over all density matrices $\sigma$ can also be done
for $\sigma$ that have rank 1. This is because the numerical range $W(X)$ is a convex set,
hence $\trace\sigma X$ and $\langle\phi, X\phi\rangle$ cover the same set.
We can thus replace (*) by
$$
\min_\phi \max_\psi \langle \psi, |X-\langle\phi,X\phi\rangle\id|_*^2\psi\rangle.
$$
A short calculation yields that this is equal to
$$
\min_\phi \max_\psi \langle \psi, |X|_*^2\psi\rangle +
|\langle\phi,X\phi\rangle - \langle\psi,X\psi\rangle|^2 - |\langle\psi,X\psi\rangle|^2,
$$
and one sees that the minimum over $\phi$ is obtained for $\phi=\psi$, and is equal to
$$
\max_\psi \langle \psi, |X|_*^2\psi\rangle - |\langle\psi,X\psi\rangle|^2,
$$
which proves that the maximum $*$-variance of $X$ over all $\rho$ is indeed obtained for $\rho$ of rank 1.
\qed

\subsection{Radius and Cartesian radius}
Because of theorem \ref{th:radiusM}, we define:
\begin{definition}
The $*$-radius of a non-normal matrix $X$ is
\end{definition}
\be
r_*(X) = \min_{y\in\C} || \,\, |X-y\id|_*\,\, ||_{\infty}.\label{eq:*radius}
\ee
where $*$ may stand for $L$, $R$ and $C$, corresponding to the use of the respective $*$-modulus.
By the theorem we've just proven, we also have the dual definition
\be
r_*(X) = \max_{\rho} (\trace[\rho|X-\trace[\rho X]\id_d|_*^2])^{1/2}. \label{eq:*radiusdual}
\ee
It is easy to see that left and right moduli yield the same value;
moreover, the Cartesian modulus yields a radius that is bounded above by the left/right radius.
\begin{theorem}\label{th:lrm}
For any matrix $X$,
$$
r_C(X) \le r_L(X) = r_R(X).
$$
\end{theorem}
\textit{Proof.}
The statement of equality of $L$ and $R$ radius follows from their definition and the fact
that $||XX^*||=||X^*X||$ for any UI norm.

Let $\hat{\rho}$ be an optimal $\rho$ in the dual expression (\ref{eq:*radiusdual}) for $r_C(X)$.
In general, $\hat{\rho}$ is not optimal for $r_L$ nor $r_R$.
Thus,
\beas
r_C^2(X) &=& \trace[\hat{\rho}|X-\trace[\hat{\rho} X]\id_d|_C^2] \\
&=& (\trace[\hat{\rho}|X-\trace[\hat{\rho} X]\id_d|_L^2] + \trace[\hat{\rho}|X-\trace[\hat{\rho} X]\id_d|_R^2])/2 \\
&\le& (r_L^2(X)+r_R^2(X))/2 \\
&=& r_L^2(X) = r_R^2(X).
\eeas
\qed

By this result, we no longer need to distinguish between $r_L(X)$ and $r_R(X)$, and we'll denote it just by $r(X)$
and call it the \textit{radius} of $X$, while we call $r_C(X)$ the \textit{Cartesian radius}.

For the proof of Theorem \ref{th:radiusM}
we needed the matrix equivalent of lemma \ref{lem:center}.
This lemma already appeared in Stampfli's paper \cite{stampfli} but was proven in a different way and only for 
the left modulus.
\begin{lemma}\label{lem:centerW}
For any matrix $X$, the value of $y\in\C$ that achieves the minimum of
$||\,|X-y\id|_*\,||$ is contained in the numerical range $W(X)$.
\end{lemma}
\textit{Proof.}
We will prove this by contradiction.
A point $z\in \C$ is in the numerical range $W(X)$ if and only if
\cite{HJII}
$$
\forall \phi\in\R: \Re(e^{i\phi}z) \le \lambda_{\max}(\Re(e^{i\phi}X))),
$$
where the real part of a matrix is defined as $\Re A = (A+A^*)/2$.

Let $y'$ be a complex number that is not in $W(X)$. Thus there exists an angle $\phi$
such that $\Re(e^{i\phi}y') > \lambda_{\max}(\Re(e^{i\phi}X)))$, strictly, or
$$
\lambda_{\max}(\Re(e^{i\phi}(X-y'\id))) < 0.
$$
We will show that this $y'$ cannot be optimal for $\min_{y} || \,\,|X-y\id|_R^2\,\,||_\infty$.

Obviously,
$|X-y\id|_R^2 = |e^{i\phi}X-e^{i\phi}y\id|_R^2$.
Thus, defining $Z=e^{i\phi}X-e^{i\phi}y'\id$ and setting $y=y'+e^{-i\phi}\epsilon$, we only need to prove that
if $\lambda_{\max}(\Re Z)<0$, then the minimum of
$\lambda_{\max}(|Z-\epsilon\id|_R^2)$ is not achieved for $\epsilon=0$.
Since this is a convex function of $\epsilon$,
it suffices to consider values of $\epsilon$ in an arbitrarily small neighbourhood of 0.

Now put $Z=-(A+iB)$, with $A$ and $B$ Hermitian.
The condition $\lambda_{\max}(\Re Z)<0$ means that $A$ should be strictly positive definite.
Does $A>0$ imply that $|| \,\,|(A+\epsilon\id)+iB|_R^2\,\,||_\infty$ is not minimal in $\epsilon=0$?
It turns out that it suffices to consider real $\epsilon$ only.
A short calculation shows
\beas
|(A+\epsilon\id)+iB|_R^2 &=&
|A+\epsilon\id|^2 - |A|^2 + |A+iB|_R^2 \\
&=& 2\epsilon(A+\epsilon\id/2) + |A+iB|_R^2.
\eeas
Since $A>0$, we can choose an $\epsilon<0$ such that we still have $A+\epsilon\id/2>0$ strictly.
Thus, there is an $\eta>0$ (given by $\lambda_{\min}(A)+\epsilon/2$)
such that $A+\epsilon\id/2>\eta\id$. Then we have
$2\epsilon(A+\epsilon\id/2)\le 2\epsilon\eta\id$.
Therefore,
\beas
\lambda_{\max}(|(A+\epsilon\id)+iB|_R^2)
&=& \lambda_{\max}(2\epsilon(A+\epsilon\id/2) + |A+iB|_R^2) \\
&\le& \lambda_{\max}(2\epsilon\eta\id + |A+iB|_R^2) \\
&=& 2\epsilon\eta + \lambda_{\max}(|A+iB|_R^2) \\
&<& \lambda_{\max}(|A+iB|_R^2).
\eeas
Thus, indeed, $\epsilon=0$ is not the minimum, as we set out to prove.

One immediately verifies that the same reasoning holds for the left modulus and the Cartesian modulus too.
\qed
\subsection{Radius compared to numerical radius}
One can now ask how these different radii $r(X)$ and $r_C(X)$
relate to the numerical range $W(X)$. While we do not know the
ultimate answer, we do know that none of the radii is the radius
of the smallest circle circumscribing $W(X)$.
The Cartesian radius of $X$ can be expressed as
$$
r_C^2(X) = \min_{z\in \C} \max_\rho \trace\rho |X-z\id|_C^2.
$$
The radius of $W(X)$ is
$$
r_W(X) := r(W(X)) = \min_{z\in \C} \max_\rho |\trace\rho (X-z\id)|.
$$
Again, this is not to be confused with the numerical radius, $w(X):=\max_\rho |\trace\rho X|$.
We therefore call $r_W$ the \textit{central numerical radius}. We have:
$$
r_W(X) = \min_{z\in\C} w(X-z\id).
$$

We now show that the central numerical radius is never bigger than the Cartesian radius.
This follows directly from:
\begin{theorem}
For all matrices $X$,
$$
w(X)\le ||\,\, |X|_C \,\,||.
$$
\end{theorem}
\textit{Proof.}
In terms of the Cartesian decomposition of $X=A+iB$,
$$
w(X) = \max_\rho |\trace\rho(A+iB)| = \max_\rho \sqrt{(\trace\rho A)^2+(\trace\rho B)^2},
$$
and
$$
||\,\, |X|_C \,\,|| =||\sqrt{A^2+B^2}|| = \max_\rho \trace\rho\sqrt{A^2+B^2}.
$$
The theorem would follow if, for all density matrices $\rho$
and Hermitian $A$ and $B$,
\be
\sqrt{(\trace\rho A)^2+(\trace\rho B)^2} \le \trace\rho\sqrt{A^2+B^2}. \label{eq:ep1}
\ee
Note first that
$|\trace\rho A| \le \trace\rho |A|$, thus we only have to prove the inequality for positive $A$ and $B$.
Indeed, let $A=A_+-A_-$ be the Jordan decomposition of $A$, then
$|\trace\rho A| = |\trace\rho A_+ - \trace\rho A_-| \le |\trace\rho_+| + |\trace\rho A_-|
=\trace\rho_+ + \trace\rho A_- = \trace\rho |A|$.

By making the substitutions $A=X^{1/2}$ and $B=Y^{1/2}$, and taking squares on both sides,
the inequality becomes
\be
(\trace\rho X^{1/2})^2+(\trace\rho Y^{1/2})^2 \le (\trace\rho\sqrt{X+Y})^2,\label{eq:ep2}
\ee
which expresses the concavity of the function $X\mapsto (\trace\rho X^{1/2})^2$ on the set of positive
matrices. It turns out that the function $X\mapsto (\trace\rho X^{1/p})^p$ is concave for all $p\ge 1$.
This can be proven by reducing the statement to
Epstein's theorem \cite{epstein}, which states that the function $X\mapsto \trace(B X^{1/p} B)^p$
is concave for all $p\ge 1$.
Taking, in particular, $B=\psi\psi^*$, with $\psi$ any normalised vector,
shows that the function $X\mapsto(\psi^* X^{1/p}\psi)^p$ is concave, and that already proves (\ref{eq:ep2})
and (\ref{eq:ep1}) for $\rho$ that have rank 1.
The validity of (\ref{eq:ep1}) for general $\rho$ then follows immediately by noting that any density
matrix $\rho$ can be written as a convex combination of rank 1 density matrices,
the left-hand side of (\ref{eq:ep1}) is convex in $\rho$, and the right-hand side is linear in $\rho$.
\qed

This easily gives:
\begin{corollary}
For all matrices $X$,
$r_W(X)\le r_C(X)$.
\end{corollary}
\textit{Proof.}
By the previous theorem, for all $z\in \C$, $w(X-z\id)\le ||\,\, |X-z\id|_C \,\,||$.
Minimising both sides over all $z\in \C$ then gives $r_W(X) \le r_C(X)$.
\qed
\subsection{Radius compared to matrix norms}
Coming back to the definition of the various radii, as given by (\ref{eq:*radius}),
one can again ask whether the infinity norm in (\ref{eq:radiusM}) has to
be replaced by the second Ky-Fan norm,
as was the case for normal matrices, to yield the best possible norm based bounds on the radii.
The answer is negative.
Instead, we have the following theorem that
gives a bound on the $L$ and $R$ radius in terms of the infinity norm, and
a bound on the Cartesian norm in terms of the Ky Fan $||.||_{(2),2}$-norm.
The reason for these different choices of norms is because these norms
turn out to be the fundamental ones for each case, from which best bounds for
every other norm can be derived.

It can be expected that non-normal matrices might allow larger radii for fixed given norm.
This is indeed the case.
The best bound for the $L$ and $R$ radius is much weaker than in the normal case,
to the point that its proof is actually trivial.
The best bound for the Cartesian radius is
stronger, and coincides with the bounds for the normal case for many norms.
To see this, compare for example corollary \ref{th:rcvp} below with theorem \ref{th:normalmain};
more precisely, the normal and non-normal bounds coincide for Schatten $p$-norms with $p\ge2$.
As could be expected, the proof is also harder.
This can be seen as an indication that the Cartesian norm is
the natural norm to use as far as radii of non-normal matrices are concerned.
\begin{theorem}\label{th:nonnormal}
For any $n\times n$ matrix $X$,
$$
r_L(X) \le ||X||_{(1)} = ||X||_{\infty},
$$
while
$$
r_C(X) \le \frac{1}{\sqrt{2}} ||X||_{(2),2}.
$$
\end{theorem}
\textit{Proof.}
The bound for $r_L$ follows immediately from the definition (\ref{eq:*radius}) by replacing the optimal $y$
by the suboptimal $y=0$.

For the $r_C$ bound, we will exploit the fact that there is a rank 1 density matrix $\hat{\rho}$ achieving
optimality in $r_C^2(X) = \max_{\rho} \trace\rho |X|_C^2-|\trace\rho X|^2$.
Let $\psi$ be the normalised vector in $\C^n$ for which $\hat{\rho}=\psi\psi^*$.
We can now construct two orthonormal bases $\{u_i\}_{i=1}^n$ and $\{v_i\}_{i=1}^n$,
with $u_1=v_1=\psi$ and all other vectors unspecified for the time being,
and express $X$ in these bases as $X=\sum_{i,j} x_{ij} u_i v_j^*$
with $x_{ij} =\langle u_i, X v_j\rangle$.
The Cartesian radius of $X$ is then given by
\beas
r_C(X)^2 &=& \frac{1}{2} \langle \psi,(XX^*+X^*X)\psi\rangle - |\langle\psi, X\psi\rangle|^2 \\
&=& \frac{1}{2}(\langle u_1,(XX^*)u_1\rangle+\langle v_1,(X^*X)v_1\rangle) - |\langle u_1, X v_1\rangle|^2 \\
&=& \frac{1}{2}\left(\sum_{j=1}^n \langle u_1,X v_j\rangle \langle v_j,X^* u_1\rangle
+\sum_{j=1}^n \langle v_1,X^* u_j\rangle \langle u_j,X v_1\rangle\right) \\
&& \mbox{ }- |\langle u_1, X v_1\rangle|^2 \\
&=& \frac{1}{2}\left(\sum_{j=1}^n |x_{1j}|^2+\sum_{j=1}^n |x_{j1}|^2\right)- |x_{11}|^2 \\
&=& \frac{1}{2}\sum_{j=2}^n(|x_{1j}|^2+|x_{j1}|^2).
\eeas

We can use the remaining degrees of freedom in the two bases for choosing their vectors in such a way that
all matrix elements $x_{1j}$ and $x_{j1}$ with $j>2$ are zero.
Then we get the simple expression $r_C(X)^2 = (|x_{12}|^2 + |x_{21}|^2)/2$.

Obviously, an upper bound on $(|x_{12}|^2 + |x_{21}|^2)/2$
is $\sum_{j=1}^2(|x_{1j}|^2+|x_{2j}|^2)/2 = ||X'||_2^2/2$,
where $X'$ is the $2\times n$ matrix consisting of the upper 2 rows of $X$ in the chosen bases.
This can be written differently: let $P$ be the $2\times n$ matrix given by
$P=e^1 u_1^*+e^2 u_2^*$, then $||X'||_2 = ||PX||_2$.
Hence, $||X'||_2^2 = \trace(PXX^*P^*)=\trace(P^*PXX^*)$.
Now note that $P^*P$ is a rank 2 partial isometry. Thus an upper bound on $||X'||_2^2$
is given by the maximum of $|\trace(AXX^*)|$ over all rank 2 partial isometries.
By Ky Fan's maximum principle, this maximum is equal to $\sigma_1(XX^*)+\sigma_2(XX^*) =
\sigma_1(X)^2+\sigma_2(X)^2$.
Therefore, $(\sigma_1(X)^2+\sigma_2(X)^2)/2$
is an upper bound on $||X'||_2^2/2$ and also on $r_C(X)^2$,
proving the second inequality of the theorem.
\qed

We obtain as a corollary:
\begin{corollary}\label{th:rcvp}
For every matrix $X$,
$$
r_L(X) \le ||X||_p,\quad p\ge1
$$
and
$$
r_C(X) \le \left\{
\begin{array}{ll}
2^{-1/p}||X||_p,& p\ge2 \\
2^{-1/2}||X||_p,& 1\le p\le2.
\end{array}
\right.
$$
These inequalities are sharp.
\end{corollary}
\textit{Proof.}
Consider first the $L$-radius.
As is well-known, $||X||_{(1)}\le ||X||_p$ for all $p\ge1$.
Equality is obtained for $X=e^{12}$.

For the $C$-radius, we have,
by Corollary \ref{cor:kyfan} with $p_0=2$,
$\frac{1}{\sqrt{2}} ||XX^*||_{(2)}^{1/2} \le ||X||_p/||F||_p = 2^{-1/p}||X||_p$
for all $p\ge 2$, so that $r_C(X)\le 2^{-1/p}||X||_p$
for all $p\ge 2$. In addition, since $||X||_p\ge||X||_2$ for $1\le p\le 2$,
we also have $r_C(X)\le 2^{-1/2}||X||_p$ for $1\le p\le 2$.

Equality for $1\le p\le 2$ is obtained for $X=e^{12}$, and for $p\ge 2$ for
$X=F$.
\qed

It would have been nice if the following had been true:
\be
r_C(X)\le ||\,|X|_C\,||_{(2)}/2,\label{eq:wrong}
\ee
since in combination with theorem \ref{th:kyfan}
this would have given
$$
r_C(X) \le |||\,|X|_C\,|||/|||F|||,
$$
and, in particular, for Schatten $p$-norms
$$
r_C(X) \le 2^{-1/p}||\,|X|_C\,||_p.
$$
In fact, for $d>2$ none of these inequalities are true. If they had been,
the Bhatia-Kittaneh inequalities (\ref{eq:bk}) would have given an alternative proof of Corollary \ref{th:rcvp}.
The fact that numerical tests showed (\ref{eq:wrong}) to hold for $d=2$ 
provided the inspiration for the proof of theorem
\ref{th:nonnormal}.
\subsection{Application to commutator bounds}
We finish by giving the promised sharp bound on the Frobenius norm of a commutator:
\begin{corollary}
For general complex matrices $X$ and $Y$, and $p\ge 1$,
$$
||\,[X,Y]\,||_2 
\le \sqrt{2} ||X||_2 ||Y||_{(2),2}
\le 2^{\max(1/2,1-1/p)}||X||_2 ||Y||_p.
$$
\end{corollary}
\textit{Proof.}
In the proof of theorem \ref{th:frobnorm} we already found that
$$
||XY-YX||_2^2 \le 4||X||_2^2 \left(\trace[\rho(Y^*Y+YY^*)/2]-|\trace[\rho Y]|^2\right).
$$
The second factor is what we coined the Cartesian variance of $Y$, $\var_{C,\rho}(Y)$, and is thus bounded
above by $r_C(Y)^2$. By theorem \ref{th:nonnormal} and corollary \ref{th:rcvp}, we find
$$
||XY-YX||_2 \le \sqrt{2} ||X||_2 ||Y||_{(2),2}
$$
and the other stated inequalities.
\qed
\begin{ack}
I am grateful for the hospitality of the following institutions where parts of this work were done:
the Banff Center, Banff, Canada and the Fields Institute, Toronto, also in Canada.
Thanks to R.F.\ Werner, S.\ Michalakis, Mary-Beth Ruskai and Michael Nathanson
for discussions, David Wenzel for sharing
his preprint, John Holbrook for the reference to Stampfli's work, and, 
finally, the staff of `Hotel Energetyk', Myczkowce, Poland,
for providing the energy.
\end{ack}

\end{document}